\documentclass[leqno,a4paper,11pt]{article}
\usepackage{amsmath}
\usepackage{amssymb}
\usepackage{amsthm}

\begin{document}

\setlength{\baselineskip}{20pt}

\newtheorem{thm}{\sc{Theorem}}[section]
\newtheorem{prop}[thm]{\sc{Proposition}}
\newtheorem{lem}[thm]{\sc{Lemma}}
\newtheorem{cor}[thm]{\sc{Corollary}}
\renewcommand{\proofname}{\sc{Proof}}

\numberwithin{equation}{section}



\renewcommand{\theenumi}{\roman{enumi}}
\renewcommand{\labelenumi}{$\mathrm{(\theenumi)}$}

\renewcommand{\refname}{\begin{center} References \end{center}}

\renewcommand{\abstractname}{}

\renewcommand{\thefootnote}{}

\renewcommand{\figurename}{\bf {Fig.}}

\setcounter{section}{-1}

\title{Quantum $b$-functions of prehomogeneous
vector spaces of commutative parabolic type}

\author{Atsushi KAMITA}

\date{}

\maketitle

\begin{abstract}
\hspace{-7mm} \textbf{\sc{Abstract.}}
We show that there exists a natural $q$-analogue of the $b$-function
for the prehomogeneous vector space of commutative parabolic type,
and calculate them explicitly in each case.
Our method of calculating the $b$-functions seems to be new even
for the original case $q=1$.
\end{abstract}

\footnotetext{2000 {\it Mathematics Subject Classification}:
Primary 17B37; Secondary 17B10, 20G05.}
\footnotetext{{\it Keywords and Phrases}: quantum groups,
semisimple Lie algebras, highest weight modules, $b$-functions.}

\section{Introduction}

Among prehomogeneous vector spaces those called of commutative parabolic
type have special features since they have additional information coming
from their realization inside simple Lie algebras.
In \cite{KMT} we constructed
a quantum analogue $A_q(V)$ of the
coordinate algebra $A(V)$ for a prehomogeneous vector space $(L,V)$
of commutative parabolic type.
If $(L, V)$ is regular, then there exists a
basic relative invariant $f \in A(V)$.
In this case a quantum analogue $f_q \in A_q(V)$ of $f$ is also
implicitly constructed in \cite{KMT}.
The aim of this paper is to give a quantum analogue of the $b$-function of $f$.

Let ${^t}f(\partial)$ be the constant coefficient differential operator
on $V$ corresponding to the relative invariant ${^t}f$
of the dual space $(L, V^{*})$.
Then the $b$-function $b(s)$ of $f$ is given by
${^t}f(\partial) f^{s+1} = b(s) f^{s}$.
See \cite{Kimura}, \cite{MRS} and \cite{H-U} for the explicit form of $b(s)$.

For $g \in A_q(V)$ we can also define a (sort of $q$-difference)
operator ${^t}g(\partial)$ by
\begin{align*}
\langle {^t}g(\partial) h, h' \rangle = \langle h, g h \rangle
\quad (h, h' \in A_q(V)),
\end{align*}
where $\langle \ , \ \rangle$ is a natural non-degenerate symmetric
bilinear form on $A_q(V)$ (see Section \ref{def of q-b-func} below).
We can show that there exists some $b_q(s) \in \mathbb{C}(q)[q^s]$
satisfying
\begin{align*}
{^t}f_q(\partial) f_q^{s+1} = b_q(s) f_q^{s} \quad (s \in \mathbb{Z}_{\ge 0}).
\end{align*}

Our main result is the following.
\begin{thm}
If we have $b(s) = \prod_{i} (s + a_i)$, then we have
\begin{align*}
b_q(s) = \prod_i q_0^{s + a_i -1 } [s + a_i]_{q_{0}}
\quad (\textrm{up to a constant multiple}),
\end{align*}
where
$q_0 = q^{2}$ $(\textrm{type $B$, $C$})$ or
$q$ $(\textrm{otherwise})$,
and $\displaystyle{[n]_t = \frac{t^{n} - t^{-n}}{t - t^{-1}}}$.
\end{thm}

We shall prove this theorem using an induction on the rank
of the corresponding simple Lie algebra.
We remark that this result was already obtained for type $A$
in Noumi-Umeda-Wakayama \cite{NUW} using a quantum analogue
of the Capelli identity.

The author expresses gratitude to Professor A.\ Gyoja and
Professor T.\ Tanisaki.

\section{Quantized enveloping algebra}

Let $\mathfrak{g}$ be a simple Lie algebra
over the complex number filed $\mathbb{C}$
with Cartan subalgebra $\mathfrak{h}$.
Let $\Delta \subset \mathfrak{h}^{*}$ be the root system and
$W \subset \mathrm{GL}(\mathfrak{h})$ the Weyl group.
For $\alpha \in \Delta$ we denote the corresponding root space by
$\mathfrak{g}_{\alpha}.$
We denote the set of positive roots by $\Delta^{+}$
and the set of simple roots by $\{\alpha_{i}\}_{i \in I_0}$,
where $I_0$ is an index set.
For $i \in I_0$ let $h_i \in \mathfrak{h}$,
$\varpi_i \in \mathfrak{h}^{*}$, $s_i \in W$ be the simple coroot,
the fundamental weight and the simple reflection corresponding to
$i$ respectively.
We denote the longest element of $W$ by $w_0$.
Let $(\ , \ ) : \mathfrak{g} \times \mathfrak{g} \rightarrow
\mathbb{C}$ be the invariant symmetric bilinear form such that
$(\alpha, \alpha) = 2$ for short roots $\alpha$.
For $i, j \in I_0$ we set
\begin{equation*}
d_i = \frac{(\alpha_i, \alpha_i)}{2},\quad
a_{i j} = \frac{2 (\alpha_i, \alpha_j)}{(\alpha_i, \alpha_i)}.
\end{equation*}
We define the antiautomorphism $x \mapsto {^t}x$ of the enveloping
algebra $U(\mathfrak{g})$ of
$\mathfrak{g}$ by
$^{t} \! x_{\alpha} = x_{-\alpha}$ and $^{t} \! h_i = h_i$,
where $\{ x_{\alpha} | \alpha \in \Delta \}$ is a Chevalley basis
of $\mathfrak{g}$.

The quantized enveloping algebra $U_q(\mathfrak{g})$
of $\mathfrak{g}$ (Drinfel'd \cite{Drinfel'd}, 
Jimbo \cite{Jimbo})
is an associative algebra over the rational function
field $\mathbb{C}(q)$ generated by the elements $\{ E_i, F_i, 
K_i^{\pm 1} \}_{i \in I_0}$ satisfying the following relations
\begin{align*}
& K_i K_j = K_j K_i, \quad
K_i K_i^{-1} = K_i^{-1} K_i = 1, \\
& K_i E_j K_i^{-1} = q_i^{a_{i j}} E_j, \quad
K_i F_j K_i^{-1} = q_i^{-a_{i j}} F_j, \\
& E_i F_j - F_j E_i = \delta_{i j}
\frac{K_i - K_i^{-1}}{q_i - q_i^{-1}}, \\
& \sum_{k = 0}^{1 - a_{i j}} (-1)^k
\begin{bmatrix}
1 - a_{i j} \\
k
\end{bmatrix}
_{q_i} E_i^{1 - a_{i j} -k} E_j E_i^k = 0 \hspace{5mm} (i \neq j), \\
& \sum_{k = 0}^{1 - a_{i j}} (-1)^k
\begin{bmatrix}
1 - a_{i j} \\
k
\end{bmatrix}
_{q_i} F_i^{1 - a_{i j} -k} F_j F_i^k = 0 \hspace{5mm} (i \neq j),
\end{align*}
where $q_i = q ^{d_i}$, and
\[
[m]_t = \frac{t^m - t^{-m}}{t - t^{-1}}, \hspace{5mm}
[m]_t \, ! = \prod_{k=1}^m [k]_t, \hspace{5mm}
\begin{bmatrix}
m \\
n
\end{bmatrix}
_t = \frac{[m]_t \, !}{[n]_t \, ! \, [m-n]_t \, !} \hspace{3mm}
(m \ge n \ge 0).
\]
For $\mu = \sum_{i \in I_0} m_i \alpha_i$ we set
$K_{\mu} = \prod_{i} K_i^{m_i}$.

We can define an algebra antiautomorphism $x \mapsto {^t}x$ of
$U_q(\mathfrak{g})$ by
\begin{align*}
&{^t} \! K_i = K_i,& &{^t} \! E_i = F_i,& &{^t} \! F_i = E_i.&
\end{align*}

We define subalgebras $U_q(\mathfrak{b}^{\pm})$,
$U_q(\mathfrak{h})$ and $U_q(\mathfrak{n}^{\pm})$
of $U_q(\mathfrak{g})$ by
\begin{align*}
& U_q(\mathfrak{b}^{+}) =
\langle K_i^{\pm 1}, E_i \, | \, i \in I_0 \rangle, &
& U_q(\mathfrak{b}^{-}) =
\langle K_i^{\pm 1}, F_i \, | \, i \in I_0 \rangle, &
& U_q(\mathfrak{h}) = \langle K_i^{\pm 1} \, | \, i \in I_0 \rangle,&\\
& U_q(\mathfrak{n}^{+}) = \langle E_i \, | \, i \in I_0 \rangle,&
& U_q(\mathfrak{n}^{-}) =
\langle F_i \, | \, i \in I_0 \rangle.& & &
\end{align*}

We set $\mathfrak{h}_{\mathbb{Z}}^{*} = \oplus_{i \in I_0}
\mathbb{Z} \varpi_{i}$.
For a $U_q(\mathfrak{h})$-module $M$ we define the weight space
$M_{\mu}$ with weight $\mu \in \mathfrak{h}_{\mathbb{Z}}^{*}$ by
\begin{align*}
M_{\mu} = \{ m \in M | K_i m = q_i^{\mu(h_i)} m \ (i \in I_0) \}.
\end{align*}

The Hopf algebra structure on $U_q(\mathfrak{g})$ is defined
as follows. The comultiplication $\Delta : U_q(\mathfrak{g}) \to
U_q(\mathfrak{g}) \otimes U_q(\mathfrak{g})$ is
the algebra homomorphism satisfying
\[
\Delta (K_i) = K_i \otimes K_i, \hspace{2mm}
\Delta (E_i) = E_i \otimes K_i^{-1} + 1 \otimes E_i, \hspace{2mm}
\Delta (F_i) = F_i \otimes 1 + K_i \otimes F_i.
\]
The counit $\epsilon : U_q(\mathfrak{g}) \to \mathbb{C}(q)$ is the
algebra homomorphism satisfying
\[
\epsilon (K_i) = 1, \hspace{3mm} \epsilon (E_i) =
\epsilon (F_i) = 0.
\]
The antipode $S : U_q(\mathfrak{g}) \to U_q(\mathfrak{g})$ is the
algebra antiautomorphism satisfying
\[
S(K_i) = K_i^{-1}, \hspace{3mm} S(E_i) = -E_i K_i, \hspace{3mm}
S(F_i) = -K_i^{-1} F_i.
\]

The adjoint action of $U_q(\mathfrak{g})$ on $U_q(\mathfrak{g})$ is
defined as follows. For $x,y \in U_q(\mathfrak{g})$ write
$\Delta (x) = \sum_k x_k^{(1)} \otimes x_k^{(2)}$ and set
$\mathrm{ad} (x) (y) = \sum_k x_k^{(1)} y S(x_k^{(2)})$.
Then $\mathrm{ad} : U_q(\mathfrak{g}) \to
\mathrm{End}_{\mathbb{C}(q)} (U_q(\mathfrak{g}))$
is an algebra homomorphism.

For $i \in I_0$ we define an algebra automorphism $T_i$ of
$U_q(\mathfrak{g})$ (see Lusztig \cite{Lusztig2}) by
\begin{align*}
& T_i(K_j) = K_j K_i^{-a_{i j}}, \\
& T_i(E_j) =
\begin{cases}
-F_i K_i & (i = j) \\
\displaystyle \sum_{k = 0}^{-a_{i j}}
(-q_i)^{-k} E_i^{(-a_{i j} - k)} E_j
E_i^{(k)} & (i \neq j),
\end{cases}
\\
& T_i(F_j) =
\begin{cases}
-K_i^{-1} E_i & (i = j) \\
\displaystyle \sum_{k = 0}^{-a_{i j}}
(-q_i)^k F_i^{(k)} F_j F_i^{(-a_{i j} -k)}
& (i \neq j),
\end{cases}
\end{align*}
where
\[
E_i^{(k)} = \frac{1}{[k]_{q_i} \, !} E_i^k, \hspace{5mm}
F_i^{(k)} = \frac{1}{[k]_{q_i} \, !} F_i^k.
\]
For $w \in W$ we choose a reduced expression
$w = s_{i_1} \dotsm s_{i_k}$, and set $T_w = T_{i_1} \dotsm T_{i_k}$.
It dose not depend on the choice of the reduced
expression by Lusztig \cite{Lusztig3}.

It is known that there exists a unique bilinear form
$(\ , \ ) : U_q(\mathfrak{b}^{-}) \times U_q(\mathfrak{b}^{+})
\rightarrow \mathbb{C}(q)$ such that for any
$x, x' \in U_q(\mathfrak{b}^{+})$,
$y, y' \in U_q(\mathfrak{b}^{-})$, and $i, j \in I_0$
\begin{align*}
&(y, x x') = (\Delta(y), x' \otimes x),&
&(y y', x) = (y \otimes y', \Delta(x)),&\\
&(K_i, K_j) = q^{-(\alpha_i, \alpha_j)},&
&(F_i, E_j) = -\delta_{i j} (q_i - q_i^{-1})^{-1},&\\
&(F_i, K_j) = 0,& &(K_i, E_j) = 0&
\end{align*}
(See Jantzen \cite{Jantzen2}, Tanisaki \cite{Tanisaki2}).

For $\mu \in \sum_{i \in I_0} \mathbb{Z}_{\ge 0} \alpha_i$
let $U_q(\mathfrak{n}^{-})_{- \mu}$ be the weight space
with weight $\mu$ relative to the adjoint action of $U_q(\mathfrak{h})$
on $U_q(\mathfrak{n}^{-})$.
For any $y \in U_q(\mathfrak{n}^{-})_{- \mu}$ and $i \in I_0$
the elements $r_i(y)$ and $r'_i(y)$ of
$U_q(\mathfrak{n}^{-})_{-(\mu - \alpha_i)}$ are defined by
\begin{align*}
&\Delta(y) \in y \otimes 1 + \sum_{i \in I_0} K_i r_i(y) \otimes
F_i +
\big(\bigoplus_{\substack{0 < \nu \le \mu\\ \nu \neq \alpha_i}}
K_{\nu} U_q(\mathfrak{n}^{-})_{-(\mu - \nu)} \otimes
U_q(\mathfrak{n}^{-})_{- \nu}\big),\\
&\Delta(y) \in K_{\mu} \otimes y + \sum_{i \in I_0} K_{\mu - \alpha_i}
F_i \otimes r'_i(y) +
\big(\bigoplus_{\substack{0 < \nu \le \mu\\ \nu \neq \alpha_i}}
K_{\mu - \nu} U_q(\mathfrak{n}^{-})_{- \nu} \otimes
U_q(\mathfrak{n}^{-})_{-(\mu - \nu)}\big).
\end{align*}

\begin{lem}\label{r_i}{\rm (see Jantzen \cite{Jantzen2})}
\begin{enumerate}
	\item We have $r_i(1) = r'_i(1) =0$ and
	$r_i(F_j) = r'_i(F_j) = \delta_{i j}$ for $j \in I_0$.
	\item We have for $y_1 \in U_q(\mathfrak{n}^-)_{- \mu_1}$ and
	$y_2 \in U_q(\mathfrak{n}^-)_{- \mu_2}$
	\begin{align*}
	r_i(y_1 y_2) = q_i^{\mu_1(h_i)} y_1 r_i(y_2) + r_i(y_1) y_2, \quad
	r'_i(y_1 y_2) = y_1 r'_i(y_2) + q_i^{\mu_2(h_i)} r'_i(y_1) y_2.
	\end{align*}
	\item We have for $x \in U_q(\mathfrak{n}^{+})$ and
	$y \in U_q(\mathfrak{n}^{-})_{- \mu}$
	\begin{align*}
	(y, E_i x) = (F_i, E_i) (r_i(y), x), \quad
	(y, x E_i) = (F_i, E_i) (r'_i(y), x).
	\end{align*}
	\item We have $\mathrm{ad}(E_i) y = (q_i - q_i^{-1})^{-1}
	(K_i r_i(y) K_i- r'_i(y))$ for $y \in U_q(\mathfrak{n}^{-})_{- \mu}$.
\end{enumerate}
\end{lem}
From Lemma \ref{r_i} (ii) we have $r_i(F_i^n) = {r'}_i(F_i^n)
= q_i^{n-1} \left[ n \right] _{q_i} F_i^{n-1}$.

\section{Commutative parabolic type}\label{commutative parabolic type}

For a subset $I$ of $I_0$ we set
\begin{align*}
&\Delta_I = \Delta \cap \sum_{i \in I} \mathbb{Z} \alpha_i,&
&\mathfrak{l}_I = \mathfrak{h} \oplus 
(\bigoplus_{\alpha \in \Delta_I} \mathfrak{g}_{\alpha}),&
&\mathfrak{n}_I^{\pm} = \bigoplus_{\alpha \in \Delta^{+} \setminus
\Delta_I} \mathfrak{g}_{\pm \alpha},&
& W_I = \langle s_i \ | \ i \in I \rangle.&
\end{align*}
Let $L_I$ be the algebraic group corresponding to $\mathfrak{l}_I$.
Assume that $\mathfrak{n}_I^{+} \neq 0$
and $[\mathfrak{n}_I^{+}, \mathfrak{n}_I^{+}] = 0$.
Then it is known that $I = I_0 \setminus \{i_0\}$ for some $i_0 \in I_0$
and $(L_I, \mathfrak{n}^{+}_{I})$ is
a prehomogeneous vector space.
Since $\mathfrak{n}_I^{-}$ is identified the dual space of
$\mathfrak{n}_I^{+}$ via the Killing form,
we have $\mathbb{C}[\mathfrak{n}_I^{+}] \simeq S(\mathfrak{n}_I^{-})
= U(\mathfrak{n}_I^{-})$.
There exists finitely many $L_I$-orbits 
$C_1$, $C_2$, \dots, $C_r$, $C_{r+1}$
on $\mathfrak{n}_I^{+}$satisfying the closure relation
$\{0 \}=C_1 \subset \overline{C_2} \subset \dots \subset
\overline{C_{r}} \subset
\overline{C_{r+1}} = \mathfrak{n}_I^{+}$.
In the remainder of this paper we denote by $r$ the number of
non-open orbits on $\mathfrak{n}_I^+$.
For $p \le r$
we set $\mathcal{I}(\overline{C_p}) =
\{ f \in \mathbb{C}[\mathfrak{n}_I^{+}] \ | \ f(\overline{C_p}) = 0 \}$.
We denote by $\mathcal{I}^{m}(\overline{C_p})$ the subspace of
$\mathcal{I}(\overline{C_p})$ consisting of homogeneous elements
with degree $m$.
It is known that
$\mathcal{I}^{p}(\overline{C_p})$ is an irreducible
$\mathfrak{l}_I$-module and
$\mathcal{I}(\overline{C_p}) = \mathbb{C}[\mathfrak{n}_I^{+}] \
\mathcal{I}^{p}(\overline{C_p})$.
Let $f_p$ be the highest weight
vector of $\mathcal{I}^{p}(\overline{C_p})$, and let $\lambda_p$ be
the weight of $f_p$.
We have the irreducible decomposition
\begin{equation*}
\mathbb{C}[\mathfrak{n}_I^{+}] = \bigoplus_{\mu \in \sum_{p=1}^r
\mathbb{Z}_{\ge 0} \lambda_p} V(\mu),
\end{equation*}
where $V(\mu)$ is an irreducible highest weight module with
highest weight $\mu$ and
$V(\lambda_p) = \mathcal{I}^{p}(\overline{C_p})$
(see Schmid \cite{Schmid} and Wachi \cite{Wachi}).

If the prehomogeneous vector space $(L_I, \mathfrak{n}_I^{+})$
is regular, there exists a one-codimensional orbit $C_r$.
Then it is known that
$\mathcal{I}^{r}(\overline{C_r}) = \mathbb{C} f_r$,
$f_r$ is the basic relative invariant of $(L_I, \mathfrak{n}_I^{+})$
and $\lambda_r = -2 \varpi_{i_0}$, where $I = I_0 \setminus \{ i_0 \}$.
The pairs $(\mathfrak{g}, i_0)$ where
$(L_I, \mathfrak{n}_I^{+})$ are regular are given by the Dynkin
diagrams of Figure \ref{regular}.
Here the white vertex corresponds to $i_0$.
\begin{figure}[tbp]
	$(A_{2n-1}, n)$ \hspace{4.5mm}
	\begin{picture}(165,15)
	\thicklines
	\put(5,0){$\stackrel{1}{\bullet}$}
	\put(10,2.5){\line(1,0){15}}
	\put(25,1.5){\ldots}
	\put(45,2.5){\line(1,0){15}}
	\put(55,0){$\stackrel{n-1}{\bullet}$}
	\put(65,2.5){\line(1,0){15}}
	\put(80,0){$\stackrel{n}{\circ}$}
	\put(85,2.5){\line(1,0){15}}
	\put(95,0){$\stackrel{n+1}{\bullet}$}
	\put(105,2.5){\line(1,0){15}}
	\put(120,1.5){\ldots}
	\put(140,2.5){\line(1,0){15}}
	\put(145,0){$\stackrel{2n-1}{\bullet}$}
	\end{picture}
	
	$(B_n,1)$ \hspace{10mm}
	\begin{picture}(192,25)
	\thicklines
	\put(5,0){$\stackrel{1}{\circ}$}
	\put(10,2.5){\line(1,0){15}}
	\put(25,0){$\stackrel{2}{\bullet}$}
	\put(30,2.5){\line(1,0){15}}
	\put(45,1.5){\ldots}
	\put(65,2.5){\line(1,0){15}}
	\put(75,0){$\stackrel{n-1}{\bullet}$}
	\put(85,0){$\Longrightarrow$}
	\put(100,0){$\stackrel{n}{\bullet}$}
	\end{picture}
	
	$(C_n,n)$ \hspace{10mm}
	\begin{picture}(192,25)
	\thicklines
	\put(5,0){$\stackrel{1}{\bullet}$}
	\put(10,2.5){\line(1,0){15}}
	\put(25,0){$\stackrel{2}{\bullet}$}
	\put(30,2.5){\line(1,0){15}}
	\put(45,1.5){\ldots}
	\put(65,2.5){\line(1,0){15}}
	\put(75,0){$\stackrel{n-1}{\bullet}$}
	\put(85,0){$\Longleftarrow$}
	\put(100,0){$\stackrel{n}{\circ}$}
	\end{picture}
	
	$(D_n,1)$ \hspace{10mm}
	\begin{picture}(192,25)
	\thicklines
	\put(5,0){$\stackrel{1}{\circ}$}
	\put(10,2.5){\line(1,0){15}}
	\put(25,0){$\stackrel{2}{\bullet}$}
	\put(30,2.5){\line(1,0){15}}
	\put(45,1.5){\ldots}
	\put(65,2.5){\line(1,0){15}}
	\put(75,0){$\stackrel{n-2}{\bullet}$}
	\put(85,2.5){\line(1,0){15}}
	\put(95,0){$\stackrel{n-1}{\bullet}$}
	\put(82.5,0){\line(0,-1){15}}
	\put(80,-20){$\bullet$}
	\put(86,-20){{\footnotesize $n$}}
	\end{picture}
	
	$(D_{2n},2n)$ \hspace{6mm}
	\begin{picture}(192,45)
	\thicklines
	\put(5,0){$\stackrel{1}{\bullet}$}
	\put(10,2.5){\line(1,0){15}}
	\put(25,0){$\stackrel{2}{\bullet}$}
	\put(30,2.5){\line(1,0){15}}
	\put(45,1.5){\ldots}
	\put(65,2.5){\line(1,0){15}}
	\put(73,0){$\stackrel{2n-2}{\bullet}$}
	\put(85,2.5){\line(1,0){15}}
	\put(90,0){$\stackrel{\ \, 2n-1}{\bullet}$}
	\put(82.5,0){\line(0,-1){15}}
	\put(80,-20){$\circ$}
	\put(86,-20){{\footnotesize $2n$}}
	\end{picture}
		
	$(E_7,1)$ \hspace{10mm}
	\begin{picture}(192,45)
	\thicklines
	\put(5,0){$\stackrel{1}{\circ}$}
	\put(10,2.5){\line(1,0){15}}
	\put(25,0){$\stackrel{2}{\bullet}$}
	\put(30,2.5){\line(1,0){15}}
	\put(45,0){$\stackrel{3}{\bullet}$}
	\put(50,2.5){\line(1,0){15}}
	\put(65,0){$\stackrel{4}{\bullet}$}
	\put(70,2.5){\line(1,0){15}}
	\put(85,0){$\stackrel{6}{\bullet}$}
	\put(90,2.5){\line(1,0){15}}
	\put(105,0){$\stackrel{7}{\bullet}$}
	\put(67.5,0){\line(0,-1){15}}
	\put(65,-20){$\bullet$}
	\put(71,-20){{\footnotesize 5}}
	\end{picture}
	\vspace{5mm}
	\caption{}
	\label{regular}
\end{figure}


Assume that $(L_I, \mathfrak{n}_I^{+})$ is regular.
For $1 \le p \le r$ we set $\gamma_p = \lambda_{p-1} - \lambda_{p}$,
where $\lambda_0 = 0$.
Then we have $\gamma_p \in \Delta^{+} \setminus \Delta_I$.
We denote the coroot of $\gamma_p$ by $h_{\gamma_p}$, and set
$\mathfrak{h}^{-} = \sum_{p=1}^{r} \mathbb{C} h_{\gamma_p}$.
We set
\begin{align*}
&\Delta_{(p)}^+ = \{ \beta \in \Delta^+ \setminus \Delta_I \ \big|
\ \beta |_{\mathfrak{h}^-} = (\gamma_j + \gamma_k)/2 \ 
\textrm{\ for some } 1 \le j \le k \le p\} \cup
\{ \gamma_1, \dots, \gamma_p \},\\
&\mathfrak{n}_{(p)}^{\pm} = \sum_{\beta \in \Delta_{(p)}^+}
\mathfrak{g}_{\pm \beta},\\
&\mathfrak{l}_{(p)} =
[\mathfrak{n}_{(p)}^{+}, \mathfrak{n}_{(p)}^{-}]
\end{align*}
(see Wachi \cite{Wachi} and Wallach \cite{Wallach}).
Note that $\alpha_{i_0} \in \Delta_{(p)}^+$ for any $p$
and $\Delta_{(r)}^+ = \Delta^+ \setminus \Delta_I$.
Then it is known that
$(L_{(p)}, \mathfrak{n}_{(p)}^{+})$ is a regular
prehomogeneous vector space of commutative parabolic type,
where $L_{(p)}$ is the subgroup of $G$ corresponding
to $\mathfrak{l}_{(p)}$.
Moreover $f_j \in \mathbb{C}[\mathfrak{n}_{(p)}]$
for $j \le p$,
and $f_p$ is a basic relative invariant of
$(L_{(p)}, \mathfrak{n}_{(p)}^{+})$.
The regular prehomogeneous vector space
$(L_{(r-1)}, \mathfrak{n}_{(r-1)}^{+})$
is described by the following.

\begin{lem}
\begin{enumerate}
	\item For $(A_{2n-1}, n)$
	we have $r = n$ and $(L_{(n-1)}, \mathfrak{n}_{(n-1)}^{+})
	\simeq (A_{2n-3}, n-1)$.
	\item For $(B_n, 1)$
	we have $r = 2$ and $(L_{(1)}, \mathfrak{n}_{(1)}^{+})
	\simeq (A_{1}, 1)$.
	\item For $(C_n, n)$ \, $(n \ge 3)$
	we have $r = n$ and $(L_{(n-1)}, \mathfrak{n}_{(n-1)}^{+})
	\simeq (C_{n-1}, n-1)$.
	\item For $(D_n, 1)$
	we have $r = 2$ and $(L_{(1)}, \mathfrak{n}_{(1)}^{+})
	\simeq (A_{1}, 1)$.
	\item For $(D_{2n}, 2n)$ \, $(n \ge 3)$
	we have $r = n$ and $(L_{(n-1)}, \mathfrak{n}_{(n-1)}^{+})
	\simeq (D_{2n-2}, 2n-2)$.
	\item For $(E_7, 1)$
	we have $r = 3$ and $(L_{(2)}, \mathfrak{n}_{(2)}^{+})
	\simeq (D_{6}, 1)$.
\end{enumerate}
\end{lem}

\section{Quantum deformations of coordinate algebras}
\label{section q-coordinate}

In this section we recall basic properties of the quantum analogue
of the coordinate algebra $\mathbb{C}[\mathfrak{n}^{+}_I]$
of $\mathfrak{n}^{+}_I$ satisfying
$[\mathfrak{n}_I^{+}, \mathfrak{n}_I^{+}] = 0$ (see \cite{KMT}).
We do not assume that $(L_I, \mathfrak{n}_I^+)$ is regular.
We take $i_0 \in I_0$ as in Section \ref{commutative parabolic type}.

We define a subalgebra $U_q(\mathfrak{l}_I)$ by
$U_q(\mathfrak{l}_I) = \langle K_i^{\pm 1}, E_j, F_j \, |
\, i \in I_0, j \in I \rangle$.
Let $w_I$ be the longest element of $W_I$,
and set
\begin{equation*}
U_q(\mathfrak{n}_I^{-}) = U_q(\mathfrak{n}^-) \cap
T_{w_I}^{-1} U_q(\mathfrak{n}^-).
\end{equation*}
We take a reduced expression $w_I w_0 =
s_{i_1} \dots s_{i_k}$ and set
\begin{align*}
&\beta_t = s_{i_1} \dotsm s_{i_{t - 1}} (\alpha_{i_t}),&
&Y_{\beta_t} = T_{i_1} \dotsm T_{i_{t - 1}} (F_{i_t})&
\end{align*}
for $t = 1, \dots , k$.
In particular $Y_{\beta_1} = F_{i_0}$.
We have $\{ \beta_t \, | \, 1 \le t \le k \} =
\Delta^{+} \setminus \Delta_I$.
The set $\{ Y_{\beta_1}^{n_1} \dotsm Y_{\beta_k}^{n_k} \, |
\, n_1, \dots, n_k \in \mathbb{Z}_{\ge 0} \}$ is a basis of
$U_q(\mathfrak{n}_I^{-})$.

\begin{prop}\label{q-coordinate}{\rm (see \cite{KMT})}
\begin{enumerate}
	\item We have $\mathrm{ad}(U_q(\mathfrak{l}_I))
	\; U_q(\mathfrak{n}_{I}^{-})
	\subset U_q(\mathfrak{n}_{I}^{-})$.
	\item The elements $Y_{\beta} \in U_q(\mathfrak{n}_I^{-})$ for
	$\beta \in \Delta^{+} \setminus \Delta_I$ do not depend on
	the choice of a reduced expression of $w_I w_0$, and they satisfy
	quadratic fundamental relations as generators of the algebra
	$U_q(\mathfrak{n}_{I}^{-})$
\end{enumerate}
\end{prop}
We regard the subalgebra $U_q(\mathfrak{n}_{I}^{-})$ of
$U_q(\mathfrak{n}^{-})$ as a quantum analogue of the coordinate
algebra $\mathbb{C}[\mathfrak{n}_I^{+}]$ of $\mathfrak{n}_I^{+}$.

Since $\mathbb{C}[\mathfrak{n}_I^+]$ is a multiplicity free
$\mathfrak{l}_I$-module, for the $L_I$-orbit $C_p$
on $\mathfrak{n}_I$ there exist unique
$U_q(\mathfrak{l}_I)$-submodules $\mathcal{I}_q(\overline{C_p})$
and $\mathcal{I}_q^{p}(\overline{C_p})$ of $U_q(\mathfrak{n}_I^-)$
satisfying
\begin{align*}
&\mathcal{I}_q(\overline{C_p}) \big|_{q=1}=
\mathcal{I}(\overline{C_p}),&
&\mathcal{I}_q^{p}(\overline{C_p}) \big|_{q=1}=
\mathcal{I}^{p}(\overline{C_p})&
\end{align*}
(see \cite{KMT}).

\begin{prop}\label{q-ideal}{\rm (see \cite{KMT})} \
$\mathcal{I}_q(\overline{C_p})=U_q(\mathfrak{n}_I^-) \;
\mathcal{I}_q^{p}(\overline{C_p})=\mathcal{I}_q^{p}(\overline{C_p})
\; U_q(\mathfrak{n}_I^-)$.
\end{prop}

Let $f_{q, p}$ be the highest weight vector of
$\mathcal{I}_q^{p}(\overline{C_p})$.
We have the irreducible decomposition
\begin{equation*}
U_q(\mathfrak{n}_I^{-}) = \bigoplus_{\mu \in \sum_{p}
\mathbb{Z}_{\ge 0} \lambda_p} V_q(\mu),
\end{equation*}
where $V_q(\mu)$ is an irreducible highest weight module with
highest weight $\mu$ and
$V_q(\lambda_p) = \mathcal{I}_q^{p}(\overline{C_p})$.
Explicit descriptions of $U_q(\mathfrak{n}_I^{-})$ and $f_{q, p}$
are given in \cite{Kamita} in the case where $\mathfrak{g}$ is classical,
and in \cite{Morita} for the exceptional cases.

Let $f$ be a weight vector of $U_q(\mathfrak{n}_I^{-})$
with the weight $- \mu$.
If $\mu \in m \alpha_{i_0} + \sum_{i \in I} \mathbb{Z}_{\ge 0}
\alpha_i$, then $f$ is an element
of $\sum_{\beta_1, \dots, \beta_m \in \Delta^{+} \setminus \Delta_I}
\mathbb{C}(q) Y_{\beta_1} \cdots Y_{\beta_m}$.
So we can define the degree of $f$ by $\deg f = m$.
In particular $\deg f_{q,p} = p$.

\section{Quantum deformations of relative invariants}

In the remainder of this paper we assume that
$(L_I, \mathfrak{n}_I^{+})$ is regular, and
$\{ i_0 \} = I_0 \setminus I$.
Then we regard the highest weight vector $f_{q,r}$ of
$\mathcal{I}_q^{r}(\overline{C_r})$ as the quantum analogue
of the basic relative invariant.
We give some properties of $f_{q,r}$ in this section.

By $\mathcal{I}_q^{r}(\overline{C_r}) = \mathbb{C}(q) f_{q,r}$
and $\lambda_r = -2 \varpi_{i_0}$, we have the following.
\begin{prop}\label{ad-act-f}
We have
\begin{align*}
&\mathrm{ad}(K_i) f_{q,r} = f_{q,r},&
&\mathrm{ad}(E_i) f_{q,r} = 0&
&\mathit{and}& &\mathrm{ad}(F_i) f_{q,r} = 0,&
\end{align*}
for any $i \in I$,
and $\mathrm{ad}(K_{i_0}) f_{q,r} = q_{i_0}^{-2}f_{q,r}$.
\end{prop}

\begin{lem}\label{r-i(Y)}
\begin{enumerate}
	\item For $i \in I$ we have $r_i(U_q(\mathfrak{n}_I^{-})) = 0$.
	\item For $\beta \in \Delta^{+} \setminus \Delta_I$ we have
	$r'_{i_0}(Y_{\beta}) = \delta_{\alpha_{i_0}, \beta}$.
\end{enumerate}
\end{lem}

\begin{proof}
(i) By Jantzen \cite{Jantzen2} we have
\begin{align*}
\{ y \in U_q(\mathfrak{n}^{-}) | r_i(y) = 0 \} =
U_q(\mathfrak{n}^{-}) \cap T_i^{-1} U_q(\mathfrak{n}^-).
\end{align*}
On the other hand we have $U_q(\mathfrak{n}_I^{-})
\subset U_q(\mathfrak{n}^{-}) \cap T_i^{-1} U_q(\mathfrak{n}^-)$
for $i \in I$.
Hence we have $r_i(U_q(\mathfrak{n}_I^{-})) = 0$ for $i \in I$.

(ii) We show the formula by induction on $\beta$.

By the definition of $r'_{i_0}$, it is clear
that $r'_{i_0}(Y_{\alpha_{i_0}}) = r'_{i_0}(F_{i_0}) = 1$.

Assume that $\beta > \alpha_{i_0}$ and the statement is proved
for any root $\beta_1$ in $\Delta^{+} \setminus \Delta_I$
satisfying $\beta_1 < \beta$.
For some $i \in I$ we can write
\begin{align*}
Y_{\beta} = c \ \mathrm{ad}(F_i) Y_{{\beta}'}
= c \ (F_i Y_{{\beta}'} - q^{-(\alpha_i, {\beta}')} Y_{{\beta}'} F_i),
\end{align*}
where ${\beta}'=\beta - \alpha_i$ and $c \in \mathbb{C}(q)$.
Hence we have
\begin{align*}
r'_{i_0}(Y_{\beta}) =
c \ (F_i \, r'_{i_0}(Y_{{\beta}'}) -
q^{(\alpha_i, \alpha_{i_0}-{\beta}')}
r'_{i_0}(Y_{{\beta}'}) F_i).
\end{align*}
If ${\beta}' = \alpha_{i_0}$, we have
$r'_{i_0}(Y_{{\beta}'}) = 1$.
If ${\beta}' \neq \alpha_{i_0}$, we have
$r'_{i_0}(Y_{{\beta}'}) = 0$ by the inductive hypothesis.
\end{proof}

\begin{prop}
The quantum analogue $f_{q,r}$ is a central element of
$U_q(\mathfrak{n}^{-})$.
\end{prop}

\begin{proof}
For $i \in I$ we have $[F_i, f_{q,r}] =
\mathrm{ad}(F_i) f_{q,r}$.
By Proposition \ref{ad-act-f} we have to show $[F_{i_0}, f_{q,r}] = 0$.

The $q$-analogue $f_{q,r}$ is a linear combination of
$Y_{\beta_1} \cdots Y_{\beta_r}$ satisfying $\sharp \{
\beta_i | \beta_i = \alpha_{i_0} \} \le 1$
(see \cite{Kamita} and \cite{Morita}).
By using Lemma \ref{r-i(Y)} it is easy to show that
$r'_{i_0}(f_{q,r}) \neq 0$ and ${r'_{i_0}}^{2}(f_{q,r}) = 0$.
Hence we have
\begin{align*}
{r'_{i_0}}^{2}(F_{i_0} f_{q,r}) = {r'_{i_0}}^{2}(f_{q,r} F_{i_0})
= (q_{i_0}^{d_{i_0}} + 1) r'_{i_0}(f_{q,r}).
\end{align*}
On the other hand there exists $c \in \mathbb{C}(q)$
such that $F_{i_0} f_{q,r} = c f_{q,r} F_{i_0}$
by Proposition \ref{q-ideal},
hence we have $(q_{i_0}^{d_{i_0}} + 1) r'_{i_0}(f_{q,r}) =
c (q_{i_0}^{d_{i_0}} + 1) r'_{i_0}(f_{q,r})$.
Therefore we obtain $c=1$.
\end{proof}

\section{$b$-functions and their quantum analogues}\label{def of q-b-func}

We recall the definition of the $b$-function.

For $h \in S(\mathfrak{n}_{I}^{+}) \simeq \mathbb{C}[\mathfrak{n}_{I}^{-}]$,
we define the constant coefficient
differential operator $h(\partial)$ by
\begin{align*}
h(\partial) \exp B(x, y) = h(y) \exp B(x, y)
\quad x \in \mathfrak{n}_{I}^{+}, y \in \mathfrak{n}_I^{-},
\end{align*}
where $B$ is the Killing form on $\mathfrak{g}$.
It is known that there exists a polynomial $b_r(s)$
called the $b$-function of the relative invariant $f_r$
such that for $s \in \mathbb{C}$
\begin{align*}
^{t} \! f_{r}(\partial) {f_{r}}^{s+1} = b_r(s) {f_{r}}^{s}.
\end{align*}
Then we have $\deg b_r = r$.
The explicit description of $b_r(s)$ is given by
\begin{align*}
&(A_{2n-1}, n)& &b_n(s) = (s + 1) (s + 2) \cdots (s + n)&\\
&(B_n, 1)& &b_2(s) = (s + 1) \left( s + \frac{2n-1}{2} \right)&\\
&(C_n, n)& &b_n(s) = (s + 1) \left( s + \frac{3}{2} \right)
\left( s + \frac{4}{2} \right)
\cdots \left( s + \frac{n + 1}{2} \right)&\\
&(D_n, 1)& &b_2(s) = (s + 1) \left( s + \frac{2n-2}{2} \right)&\\
&(D_{2n}, 2n)& &b_n(s) = (s + 1) (s + 3) \cdots (s + 2n - 1)&\\
&(E_7, 1)& &b_3(s) = (s + 1) (s + 5) (s + 9)&
\end{align*}
(see \cite{Kimura}, \cite{MRS} and \cite{H-U}).

We define a symmetric non-degenerate bilinear form
$\langle \ , \ \rangle$
on $S(\mathfrak{n}_{I}^{-}) \simeq
\mathbb{C}[\mathfrak{n}_{I}^{+}]$ by
$\langle f, g \rangle = (^{t} \! g(\partial)f) (0)$.

\begin{lem}{\rm (see Wachi \cite{Wachi})}
For $f, \ g, \ h \in S(\mathfrak{n}_{I}^{-})
\simeq \mathbb{C}[\mathfrak{n}_{I}^{+}]$ we have
\begin{enumerate}
	\item$\langle \mathrm{ad}(u)f, g \rangle =
	\langle f, \mathrm{ad}(^{t} \! u)g \rangle$
	for $u \in U(\mathfrak{l}_I)$,
	\item $\langle f, g h \rangle = \langle ^{t} \! g(\partial)f, h
	\rangle$. 
\end{enumerate}
\end{lem}

We have
for $\beta, {\beta}' \in \Delta^{+} \setminus \Delta_I$
\begin{align*}
\langle x_{-\beta}, x_{-{\beta}'}\rangle
= \delta_{\beta, {\beta}'} \frac{2}{(\beta, \beta)}.
\end{align*}
The comultiplication $\Delta$ of $U(\mathfrak{g})$ is defined by
$\Delta(x) = x \otimes 1 + 1 \otimes x$ for $x \in \mathfrak{g}$.
We define the algebra homomorphism $\widetilde{\Delta}$ by
$\widetilde{\Delta}(x) = \tau \Delta({^t}x)$,
where $x \in U(\mathfrak{g})$ and $\tau(y_1 \otimes y_2) =
{^t}y_1 \otimes \; {^t}y_2$.
Since ${^t}x_{-\beta}(\partial)(f g) =
{^t}x_{-\beta}(\partial) (f) \> g + f \> {^t}x_{-\beta}(\partial)(g)$,
we have
\begin{align*}
\langle f g, h \rangle =
\langle f \otimes g, \widetilde{\Delta}(h) \rangle.
\end{align*}

We shall define the $q$-analogue of the differential
operator ${^t}f(\partial)$ using the $q$-analogue of
$\langle \ , \ \rangle$.

We define the bilinear form $\langle \ , \ \rangle$ on
$U_q(\mathfrak{n}_{I}^{-})$ by
\begin{align*}
\langle f, g \rangle = (q^{-1} - q)^{\deg f} (f, {^t} \! g),
\end{align*}
for the weight vectors $f, g$ of $U_q(\mathfrak{n}_I^{-})$.
It is easy to show that this bilinear form $\langle \ , \ \rangle$
is symmetric. We have the following.

\begin{prop}\label{property of bilinear form}
Let $f, g, h \in U_q(\mathfrak{n}_I^{-})$.
\begin{enumerate}
	\item $\langle f g, h \rangle =
	\langle f \otimes g, \widetilde{\Delta}(h) \rangle$,
	where $\widetilde{\Delta}(h) = \tau \Delta({^t}h)$ and
	$\tau(h_1 \otimes h_2) = {^t}h_1 \otimes \; {^t}h_2$.
	\item For $u \in U_q(\mathfrak{l}_I)$ we have
	\begin{align*}
	\langle \mathrm{ad}(u)f, g \rangle =
	\langle f, \mathrm{ad}({^t}u)g \rangle.
	\end{align*}
	\item The bilinear form $\langle \ , \ \rangle$ is non-degenerate.
\end{enumerate}
\end{prop}

\begin{proof}
(i) It is clear from the definition.

(ii) It is sufficient to show that the statement holds
for the weight vectors $f, g$ and the canonical generator $u$
of $U_q(\mathfrak{l}_I)$.
If $u = K_i$ for $i \in I_0$, then the assertion is obvious.

Let $u = E_i$ for $i \in I$.
By Lemma \ref{r_i} and Lemma \ref{r-i(Y)} we have
\begin{align*}
(\mathrm{ad}(E_i) f, {^t} \! g) = (q_i^{-1} - q_i)^{-1}
(r'_i(f), {^t} \! g) = (f, {^t} \! g E_i).
\end{align*}
On the other hand we have
\begin{align*}
(f, {^t} \! (\mathrm{ad}(F_i) g)) =
(f, {^t} \! g E_i - q_i^{-\mu(h_i)} E_i {^t} \! g),
\end{align*}
where $-\mu$ is the weight of $g$.
Since $(U_q(\mathfrak{n}_I^{-}), E_i U_q(\mathfrak{n}^{+})) = 0$
by Lemma \ref{r_i} and Lemma \ref{r-i(Y)}, we have
$(\mathrm{ad}(E_i) f, {^t} \! g) = (f, {^t} \! (\mathrm{ad}(F_i) g))$.
We have $\deg f = \deg \left(\mathrm{ad}(F_i) f \right)$, and
hence the statement for
$u = E_i$ holds.
By the symmetry of $\langle \ , \ \rangle$ it also holds
for $u = F_i$.

(iii) We take the reduced expression $w_{0} = s_{i_1}
\cdots s_{i_k} s_{i_{k+1}} \cdots s_{i_l}$ such that
$w_I w_0 = s_{i_1} \cdots s_{i_k}$.
We define $Y_{\beta_{j}}$ as in Section \ref{section q-coordinate}.
Then $\{ Y_{\beta_1}^{n_1} \cdots Y_{\beta_k}^{n_k} Y_{\beta_{k+1}}^{n_{k+1}}
\cdots Y_{\beta_{l}}^{n_l} \}$ is a basis of $U_q(\mathfrak{n}^{-})$,
and for $j > k$ we have
$Y_{\beta_j} \in U_q(\mathfrak{n}^{-}) \cap U_q(\mathfrak{l}_I)$.
Hence we have $U_q(\mathfrak{n}^{-}) = U_q(\mathfrak{n}_I^-)
+ \sum_{i \in I} U_q(\mathfrak{n}^{-}) F_i$.

Since ${^t} U_q(\mathfrak{n}^{-}) = U_q(\mathfrak{n}^{+})$,
we have $U_q(\mathfrak{n}^+) = {^t} U_q(\mathfrak{n}_I^-) +
\sum_{i \in I} E_i U_q(\mathfrak{n}^+)$.
Moreover, we have
$(U_q(\mathfrak{n}_I^{-}), E_i U_q(\mathfrak{n}^+)) = 0$
for $i \in I$.
Hence if $\langle f, g \rangle = 0$ for any
$g \in U_q(\mathfrak{n}_I^{-})$, then $(f, u) = 0$
for any $u \in U_q(\mathfrak{n}^{+})$.
Thus the assertion follows from the non-degeneracy of $(\ , \ )$.
\end{proof}

\begin{prop}\label{bilinear of Y}
For $\beta, {\beta}' \in \Delta^{+} \setminus \Delta_I$
we have
\begin{align*}
\langle Y_{\beta}, Y_{{\beta}'} \rangle =
\delta_{\beta, {\beta}'}
\left[ \frac{(\beta, \beta)}{2} \right] _{q}^{-1}.
\end{align*}
\end{prop}

\begin{proof}
By the definition it is clear that $\langle Y_{\beta}, Y_{{\beta}'} \rangle =
0$ if $\beta \neq {\beta}'$.
In the case where $\beta = {\beta}'$
we shall show the statement by the induction on $\beta$.

Since $Y_{\alpha_{i_0}} = F_{i_0}$, we obtain
$\langle Y_{\alpha_{i_0}}, Y_{\alpha_{i_0}} \rangle =
\left[ \frac{(\alpha_{i_0}, \alpha_{i_0})}{2} \right] _{q}^{-1}$.

Assume that $\beta > \alpha_{i_0}$ and the statement holds for
any root $\beta_1$ in $\Delta^{+} \setminus \Delta_I$ satisfying
$\beta_1 < \beta$.
Then there exists a root $\gamma \; (< \beta)$ in
$\Delta^{+} \setminus \Delta_I$ such that
\begin{align*}
&Y_{\beta} = c_{\gamma, \beta} \; \mathrm{ad}(F_i) Y_{\gamma},&
&Y_{\gamma} = {c'}_{\gamma, \beta} \; \mathrm{ad}(E_i) Y_{\beta},&
\end{align*}
where $i \in I$ satisfying $\beta = \gamma + \alpha_i$ and
$c_{\gamma, \beta}, {c'}_{\gamma, \beta} \in {\mathbb{C}(q)}^{*}$.
We denote by $R$ the set of the pairs $\{ \gamma, \beta \}$ as above.
By Proposition \ref{property of bilinear form}
we have for $\{ \gamma, \beta \} \in R$
\begin{align*}
\langle Y_{\beta}, Y_{\beta} \rangle =
\langle Y_{\beta}, c_{\gamma, \beta} \; \mathrm{ad}(F_i) Y_{\gamma} \rangle
= c_{\gamma, \beta} \; \langle \mathrm{ad}(E_i) Y_{\beta}, Y_{\gamma} \rangle
= \frac{c_{\gamma, \beta}}{{c'}_{\gamma, \beta}}
\langle Y_{\gamma}, Y_{\gamma} \rangle =
\frac{c_{\gamma, \beta}}{{c'}_{\gamma, \beta}}
\left[ \frac{(\gamma, \gamma)}{2} \right] _{q}^{-1}.
\end{align*}
On the other hand we have for $\{ \gamma, \beta \} \in R$
\begin{align*}
&c_{\gamma, \beta} = {c'}_{\gamma, \beta} = 1&
&{\textrm {if }} (\beta, \beta) = (\gamma, \gamma),&\\
&c_{\gamma, \beta} = (q + q^{-1})^{-1}, \ {c'}_{\gamma, \beta} = 1&
&{\textrm {if }} 4 = (\beta, \beta) > (\gamma, \gamma) = 2,&\\
&c_{\gamma, \beta} = 1, \ {c'}_{\gamma, \beta} = (q + q^{-1})^{-1}&
&{\textrm {if }} 2 = (\beta, \beta) < (\gamma, \gamma) = 4&\\
\end{align*}
(see \cite{Kamita} and \cite{Morita}).
Hence we obtain $\langle Y_{\beta}, Y_{\beta} \rangle =
\left[ \frac{(\beta, \beta)}{2} \right] _{q}^{-1}$.
\end{proof}

By Proposition \ref{property of bilinear form} and
\ref{bilinear of Y}
we can regard $\langle \ , \ \rangle$ on
$U_q(\mathfrak{n}_I^{-})$ as the $q$-analogue of
$\langle \ , \ \rangle$ on
$S(\mathfrak{n}_I^{-}) \simeq \mathbb{C}[\mathfrak{n}_I^{+}]$.

\begin{prop}\label{explicit of q-deff}
\begin{enumerate}
	\item For any $g \in U_q(\mathfrak{n}_I^{-})$ there exists a unique
	${^t} g(\partial) \in \mathrm{End}_{\mathbb{C}(q)}
	(U_q(\mathfrak{n}_I^{-}))$
	such that $\langle {^t} \! g(\partial) f, h \rangle =
	\langle f, g h \rangle$ for any $f, h \in U_q(\mathfrak{n}_I^{-})$.
	In particular we have
	\begin{align*}
	{^t} Y_{\alpha_{i_0}}(\partial) = [d_{i_0}]_q^{-1} r'_{i_0},
	\end{align*}
	and for $\beta > \alpha_{i_0}$
	\begin{align*}
	{^t} Y_{\beta}(\partial) = c_{{\beta}', \beta}
	({^t} Y_{{\beta}'}(\partial)
	\mathrm{ad}(E_i) - q_i^{-{\beta}'(h_i)} \mathrm{ad}(E_i) 
	{^t} Y_{{\beta}'}(\partial)),
	\end{align*}
	where $Y_{\beta} = c_{{\beta}', \beta} \; \mathrm{ad}(F_i) Y_{{\beta}'}$.
	\item For $f \in U_q(\mathfrak{n}_I^{-})_{-\mu}$ and
	$g \in U_q(\mathfrak{n}_I^{-})_{-\nu}$ we have 
	${^t} \! g(\partial) f \in U_q(\mathfrak{n}_I^{-})_{-(\mu-\nu)}$.
\end{enumerate}
\end{prop}

\begin{proof}
(i) The uniqueness follows from
the non-degeneracy of $\langle \ , \ \rangle$.
If there exist ${^t} g(\partial)$ and
${^t} g'(\partial)$, then we have ${^t} (g g')(\partial) =
{^t} g'(\partial) \, {^t} g(\partial)$.
Therefore we have only to show the existence of ${^t} Y_{\beta}(\partial)$
for any $\beta \in \Delta^{+} \setminus \Delta_I$.
By Lemma \ref{r_i} we have
${^t} Y_{\alpha_{i_0}}(\partial) = [d_{i_0}]_q^{-1} r'_{i_0}$.
Let $\beta > \alpha_{i_0}$. Then there exists a root
${\beta}' (< \beta)$ such that
$Y_{\beta} = c_{{\beta}', \beta} \; \mathrm{ad}(F_i) Y_{{\beta}'}$
$(c_{{\beta}', \beta} \in \mathbb{C}(q))$.
By Proposition \ref{property of bilinear form}
we can show that ${^t} Y_{\beta}(\partial) = c_{{\beta}', \beta}
({^t} Y_{{\beta}'}(\partial)
\mathrm{ad}(E_i) - q_i^{-{\beta}'(h_i)} \mathrm{ad}(E_i) 
{^t} Y_{{\beta}'}(\partial))$ easily.

(ii) The assertion follows from (i).
\end{proof}

\begin{lem}\label{ad and q-deff}
For $i \in I$
$\mathrm{ad}(E_i) {^t} \! f_{q,r}(\partial) =
{^t} \! f_{q,r}(\partial) \mathrm{ad}(E_i)$ and
$\mathrm{ad}(F_i) {^t} \! f_{q,r}(\partial) =
{^t} \! f_{q,r}(\partial) \mathrm{ad}(F_i)$.
\end{lem}

\begin{proof}
Let $y_1$, $y_2 \in U_q(\mathfrak{n}_I^{-})$.
Since $\mathrm{ad}(F_i) f_{q,r} = 0$ for $i \in I$,
we have $\mathrm{ad}(F_i) (f_{q,r} y_2) = f_{q,r} \ \mathrm{ad}(F_i) y_2$.
Hence we obtain
\begin{align*}
\langle \mathrm{ad}(E_i) {^t} \! f_{q,r}(\partial) (y_1), y_2 \rangle
= \langle y_1, f_{q,r} \ \mathrm{ad}(F_i) y_2 \rangle
= \langle y_1, \mathrm{ad}(F_i) (f_{q,r} y_2) \rangle
= \langle {^t} \! f_{q,r}(\partial) \mathrm{ad}(E_i) (y_1), y_2 \rangle.
\end{align*}
Similarly we obtain $\mathrm{ad}(F_i) {^t} \! f_{q,r}(\partial) =
{^t} \! f_{q,r}(\partial) \mathrm{ad}(F_i)$
\end{proof}

By Proposition \ref{explicit of q-deff} and Lemma \ref{ad and q-deff}
the element
${^t} \! f_{q,r}(\partial)(f_{q,r}^{s + 1})$
$(s \in \mathbb{Z}_{\ge 0})$ is the highest weight vector with
highest weight $s \lambda_r = -2 s \varpi_{i_0}$.
Since $U_q(\mathfrak{n}_I^{-})$ is a multiplicity free
$U_q(\mathfrak{l}_I)$-module, there exists
$\tilde{b}_{q,r,s} \in \mathbb{C}(q)$ such that
\begin{align*}
{^t} \! f_{q,r}(\partial)(f_{q,r}^{s + 1}) =
\tilde{b}_{q,r,s} f_{q,r}^{s}.
\end{align*}

\begin{prop}
There exists a polynomial $\tilde{b}_{q,r}(t) \in \mathbb{C}(q)[t]$
such that $\tilde{b}_{q,r,s} = \tilde{b}_{q,r}(q_{i_0}^s)$ for any
$s \in \mathbb{Z}_{\ge 0}$.
\end{prop}

\begin{proof}
Let $\psi = \psi_1 \cdots \psi_m$, where $\psi_j = r'_{i_0}$ or
$\mathrm{ad}(E_i)$ for some $i \in I$.
Set $n = n(\psi) = \sharp \{j | \psi_j = r'_{i_0} \}$.
For $k \in \mathbb{Z}_{\ge 0}$ and
$y \in U_q(\mathfrak{n}_I^{-})_{- \mu}$ we have
\begin{align*}
r'_{i_0}(f_{q,r}^{k} y)
= q_{i_0}^{k-1 + \mu(h_{i_0})} \left[ k \right] _{q_{i_0}}
f_{q,r}^{k-1} \;
r'_{i_0}(f_{q,r}) y + f_{q,r}^{k} \; r'_{i_0}(y)
\end{align*}
by the induction on $k$.
Note that $q_{i_0}^{k-1 + \mu(h_{i_0})} \left[ k \right] _{q_{i_0}}
= (q_{i_0} - q_{i_0}^{-1})^{-1} q_{i_0}^{\mu(h_{i_0})-1}
((q_{i_0}^{k})^{2} - 1)$.
Moreover
$\mathrm{ad}(E_i)(f_{q,r}^{\; k} y) =
f_{q,r}^{\; k} \; \mathrm{ad}(E_i) y$ for $i \in I$.
Hence we have
\begin{align*}
\psi(f_{q,r}^{s + 1}) =
\sum_{p = 1}^n c_{p}(q_{i_0}^{s}) f_{q,r}^{s + 1 - p} y_p,
\end{align*}
where $c_{p} \in \mathbb{C}(q)[t]$ and $y_p \in U_q(\mathfrak{n}_I^{-})$
does not depend on $s$.

By Proposition \ref{explicit of q-deff}
${^t} \! f_{q,r}(\partial)$ is a linear combination of such
$\psi$ satisfying $n(\psi) = r$.
The assertion is proved.
\end{proof}

We set $b_{q,r}(s) = \tilde{b}_{q,r}(q_{i_0}^s)$ for simplicity.
By definition we have
\begin{align*}
\langle f_{q,r}^{s+1}, f_{q,r}^{s+1} \rangle = b_{q,r}(s)
b_{q,r}(s-1) \cdots b_{q,r}(0).
\end{align*}

\section{Explicit forms of quantum $b$-functions}

Our main result is the following.

\begin{thm}\label{explicit form}
Let $b_r(s) = \prod_{i=1}^{r} (s + a_i)$ be a $b$-function
of the basic relative invariant of the regular
prehomogeneous vector space $(L_I, \mathfrak{n}_I^{+})$.
Then the quantum analogue $b_{q,r}(s)$ of $b_{r}(s)$
is given by
\begin{align*}
b_{q,r}(s) = \prod_{i=1}^{r} q_{i_0}^{s + a_i -1}
\left[ s + a_i \right] _{q_{i_0}} \
({\textit {up to a constant multiple}}),
\end{align*}
where $\{ i_0 \} = I_0 \setminus I$.
\end{thm}

We prove this theorem by calculating $b_{q,r}(s)$ in each case.

For $p = 1, \dots, r$ we define
$\Delta_{(p)}^{+}$, $L_{(p)}$ and $\mathfrak{n}_{(p)}^{\pm}$
as in Section \ref{commutative parabolic type}.
We define the subalgebra $U_q(\mathfrak{n}_{(p)}^{-})$
of $U_q(\mathfrak{n}_I^{-})$ by
\begin{align*}
U_q(\mathfrak{n}_{(p)}^{-}) = \langle Y_{\beta} \, | \, \beta \in
\Delta_{(p)}^{+} \rangle.
\end{align*}
Then $U_q(\mathfrak{n}_{(p)}^{-})$ is a $q$-analogue
of $\mathbb{C}[\mathfrak{n}_{(p)}^{+}]$, and
$f_{q,p} \in U_q(\mathfrak{n}_{(p)}^{-})$ is a $q$-analogue
of basic relative invariant $f_p$ of the regular
prehomogeneous vector space $(L_{(p)}, \mathfrak{n}_{(p)}^{+})$.
We denote by $b_{q,p}(s)$
the $q$-analogue of the $b$-function of $f_p$.

The regular prehomogeneous vector space
$(L_{(1)}, \mathfrak{n}_{(1)}^{+})$ is of type
$(A_1, 1)$, and we have $U_q(\mathfrak{n}_{(1)}^{-}) =
\langle F_{i_0} \rangle$, $f_{q,1}=  c F_{i_0}$
$(c \in \mathbb{C}(q)^{*})$.
Since ${r'}_{i_0}(F_{i_0}^{s+1}) = q_{i_0}^{s}
\left[ s+1 \right] _{q_{i_0}} F_{i_0}^{s}$, we obtain
\begin{align*}
b_{q,1}(s) = c^2 \left[ d_{i_0} \right]_q^{-1} q_{i_0}^{s}
\left[ s + 1 \right]_{q_{i_0}}.
\end{align*}

If we determine $a_p(s) \in \mathbb{C}(q)$ by
\begin{align*}
\langle f_{q,p}^{s}, f_{q,p}^{s} \rangle = a_p(s)
\langle f_{q,p-1}^{s}, f_{q,p-1}^{s} \rangle,
\end{align*}
then we have $b_{q,p}(s) = \displaystyle{
\frac{a_p(s+1)}{a_p(s)} b_{q,p-1}(s)}$.
Therefore we can inductively obtain the explicit form of $b_{q,r}$.

The next lemma is useful for the calculation of $a_p(s)$.
\begin{lem}\label{lem of calulate-1}
\begin{enumerate}
	\item For $\beta \in \Delta^+ \setminus \Delta_I$ we have
	\begin{align*}
	{^t}Y_{\beta}(\partial)(f_{q,r}^n y) = {^t}Y_{\beta}(\partial)(f_{q,r}^n)
	\mathrm{ad}(K_{\beta}^{-1}) y +
	f_{q,r}^n \ {^t} \! Y_{\beta}(\partial) y
	\quad (y \in U_q(\mathfrak{n}_I^{-})).
	\end{align*}
	\item ${^t}Y_{\beta}(\partial)(f_{q,r}^{n}) = 
	q_{i_0}^{n-1} [n]_{q_{i_0}} f_{q,r}^{n-1}
	\ {^t}Y_{\beta}(\partial)(f_{q,r})$.
	\item For $\beta \in \Delta_{(p)}^+ \setminus \Delta_{(p-1)}^+$
	we have ${^t}Y_{\beta}(\partial)(f_{q,p-1}^n) = 0$.
\end{enumerate}
\end{lem}

\begin{proof}
(i) This is proved easily by the induction on $\beta$.
Note that $\mathrm{ad}(E_i) (f_{q,r}) = 0$ for $i \in I$.

(ii) Since $f_{q,r}$ is a central element of
$U_q(\mathfrak{n}_I^{-})$, this follows from (i).

(iii) Let $\beta \in \Delta_{(p)}^+ \setminus \Delta_{(p-1)}^+$.
Then there exists some $j \in I$ such that
$\beta \in \mathbb{Z}_{>0} \alpha_j +
\sum_{i \neq j} \mathbb{Z}_{\ge 0} \alpha_i$ and
$\gamma \in \sum_{i \neq j} \mathbb{Z}_{\ge 0} \alpha_i$ for any
$\gamma \in \Delta_{(p-1)}^+$.
Hence we have $U_q(\mathfrak{n}_{(p-1)}^{-})_{-(\lambda_{p-1}-\beta)}
= \{ 0 \}$, and the statement follows.
\end{proof}

Let us give $a_p(s)$ in each case.

Let $(L_I, \mathfrak{n}_I^{+})$ be the regular prehomogeneous vector
space of type $(A_{2n-1}, n)$.
Then the number of non-open orbits $r$ is equal to $n$, and
$d_{i_0} = d_n = 1$.
Here we label the vertices of the Dynkin diagram as in Figure \ref{regular}.

Let $1 \le i, j \le n$.
We set $\beta_{i j} = \alpha_{n-i+1} + \alpha_{n-i+2} + \cdots +
\alpha_{n+j-1}$, and $Y_{i j} = Y_{\beta_{i j}}$.
For two sequences $1 \le i_1 < \dots < i_p \le n$, 
$1 \le j_1 < \dots < j_p \le n$ we set
\begin{align*}
(i_1, \dots, i_p | j_1, \dots, j_p) = 
\sum_{\sigma \in S_p} (-q)^{l(\sigma)} Y_{i_1,j_{\sigma(1)}}
\cdots Y_{i_p, j_{\sigma(p)}}.
\end{align*}
Then we have $f_{q,p} = (1, \dots, p | 1, \dots, p)$
(see \cite{Kamita}).
It is easy to show the following formula.
\begin{align}
f_{q,p} = 
\sum_{k=1}^{p} (-q^{-1})^{p-k} Y_{p,k}
(1, \dots, p-1 | 1, \dots, \check{k}, \dots, p) \label{equat}.
\end{align}
Note that $\beta_{p, k} \in \Delta_{(p)}^{+}
\setminus \Delta_{(p-1)}^{+}$
and $(1, \dots, p-1 | 1, \dots, \check{k}, \dots, p)
= \mathrm{ad}(F_{n+k} \cdots F_{n+p-1}) f_{q,p-1}$
in \eqref{equat}.

Since $\beta_{p,i} \in \Delta_{(p)}^+ \setminus \Delta_{(p-1)}^+$
for $1 \le i \le p$,
by Lemma \ref{lem of calulate-1} we have 
\begin{align*}
{^t}Y_{p, i}(\partial)(f_{q,p}^{s_1} f_{q,p-1}^{s_2}) =
q^{s_1 -1} [s_1]_{q} f_{q,p}^{s_1 - 1} \
{^t}Y_{p, i}(\partial) (f_{q,p}) \ \mathrm{ad}(K_{\beta_{p, i}}^{-1})
(f_{q,p-1}^{s_2}).
\end{align*}
On the other hand we have the following.
\begin{lem}
\begin{align*}
{^t}Y_{i, j}(\partial) f_{q,p}
= (-q)^{i+j-2} (1, \dots, \check{i}, \dots, p | 1, \dots, \check{j},
\dots, p) \hspace{5mm} (1 \le i, j \le p).
\end{align*}
\end{lem}
\begin{proof}
By Proposition \ref{explicit of q-deff}, 
the statement is proved
by the induction on $i, j$ easily.
\end{proof}

Since $\mathrm{ad}(K_{\beta_{p,i}}^{-1})
(f_{q,p-1}) = q f_{q,p-1}$ if $i \le p-1$ and
$f_{q,p-1}$ if $i = p$, we have
\begin{align*}
&\langle f_{q,p}^{s_1} f_{q, p-1}^{s_2}, f_{q,p}^{s_1} f_{q, p-1}^{s_2}
\rangle& &=& &\sum_{i=1}^{p} (-q^{-1})^{p-i} \langle 
{^t} Y_{p, i}(\partial) (f_{q,p}^{s_1} f_{q, p-1}^{s_2}), 
g_i f_{q,p}^{s_1 - 1}
f_{q, p-1}^{s_2} \rangle&\\
& & &=& &\sum_{i=1}^{p-1} (-q)^{2i-2} q^{s_1 + s_2 - 1} [s_1]_q
\langle f_{q,p}^{s_1 -1} g_i
f_{q, p-1}^{s_2}, f_{q,p}^{s_1 - 1}
g_i f_{q, p-1}^{s_2} \rangle&\\
& & & & &+(-q)^{2p-2} q^{s_1 -1} [s_1]_q
\langle f_{q,p}^{s_1 - 1} f_{q,p-1}^{s_2 + 1},
f_{q,p}^{s_1 - 1} f_{q,p-1}^{s_2 + 1} \rangle,&
\end{align*}
where $g_i = (1, \dots, p - 1 | 1, \dots, \check{i}, \dots, p)$.

Now we have for $1 \le i \le p-1$
\begin{align*}
&g_i = \mathrm{ad}(F_{n+i}) g_{i+1},&
&g_{i+1} = \mathrm{ad}(E_{n+i}) g_i,&\\
&\mathrm{ad}(E_{n+i})f_{q,p} = f_{q.p-1} = 0,&
&\mathrm{ad}(F_{n+i})f_{q,p} = 0,&
&\mathrm{ad}(F_{n+i})f_{q, p-1} = \delta_{i, p-1} g_{p-1}&
\end{align*}
(see \cite{Kamita}).
Therefore we have
\begin{align*}
\mathrm{ad}(E_{n+p-1} \cdots E_{n+i+1} E_{n+i})
(f_{q,p}^{s_1 -1} g_i f_{q, p-1}^{s_2}) =
\mathrm{ad}(E_{n+p-1} \cdots E_{n+i+1})
(f_{q,p}^{s_1 -1} g_{i+1} f_{q, p-1}^{s_2})\\
= \dots =
\mathrm{ad}(E_{n+p-1}) (f_{q,p}^{s_1 -1} g_{p-1} f_{q, p-1}^{s_2})
= q^{-s_2} f_{q,p}^{s_1 -1} f_{q, p-1}^{s_2 + 1},
\end{align*}
and
\begin{align*}
f_{q,p}^{s_1 - 1} g_i f_{q, p-1}^{s_2} =
\mathrm{ad}(F_{n+i})(f_{q,p}^{s_1 - 1} g_{i+1} f_{q, p-1}^{s_2}) =
\dots = \mathrm{ad}(F_{n+i} \cdots F_{n+p-2})
(f_{q,p}^{s_1 - 1} g_{p-1} f_{q, p-1}^{s_2}).
\end{align*}
Here we have $g_{p-1} f_{q,p-1} = q^{-1} f_{q, p-1} g_{p-1}$, and
hence
\begin{align*}
f_{q,p}^{s_1 - 1} g_i f_{q, p-1}^{s_2} =
q^{-s_2} [s_2 + 1]_q^{-1}
\mathrm{ad}(F_{n+i} \cdots F_{n+p-2} F_{n+p-1})
(f_{q,p}^{s_1 - 1} f_{q, p-1}^{s_2 + 1}).
\end{align*}
By Proposition \ref{property of bilinear form}
we obtain
\begin{align*}
\langle f_{q,p}^{s_1} f_{q, p-1}^{s_2}, f_{q,p}^{s_1} f_{q, p-1}^{s_2}
\rangle = q^{s_1 -1} [s_1]_q
(q^{-s_2}[s_2 + 1]_q^{-1} \sum_{i=1}^{p-1} q^{2i-2} + q^{2p-2})
\langle f_{q,p}^{s_1 - 1} f_{q, p-1}^{s_2 + 1},
f_{q,p}^{s_1 - 1} f_{q, p-1}^{s_2 + 1} \rangle\\
= q^{p + s_1 - 2} [s_1]_q [p + s_2]_q [s_2 + 1]_q^{-1}
\langle f_{q,p}^{s_1 - 1} f_{q, p-1}^{s_2 + 1},
f_{q,p}^{s_1 - 1} f_{q, p-1}^{s_2 + 1} \rangle.
\end{align*}
From this formula we have the following.
\begin{prop}
Let $(L_I, \mathfrak{n}_I^{-})$ be a regular prehomogeneous vector
space of type $(A_{2n-1}, n)$.
We have
\begin{align*}
a_p(s) = q^{\frac{s (s + 2p - 3)}{2}} \prod_{i=1}^{s}
[i+p-1]_q.
\end{align*}
In particular $b_{q,p}(s) = q^{s + p - 1} [s + p]_q \;b_{q,p-1}(s)$, and
we have the quantum $b$-function
\begin{align*}
b_{q,n}(s) = \prod_{p=1}^{n} q^{s + p -1} [s + p]_q.
\end{align*}
\end{prop}

\vspace{5mm}

Next, we assume that $(L_I, \mathfrak{n}_I^{+})$ is
regular of type $(D_{2n}, 2n)$.
We label the vertices of the Dynkin diagram as in
Figure \ref{regular}, then $d_{i_0} = d_{2n} = 1$.
There exist $n$ non-open orbits on $\mathfrak{n}_I^{+}$.
Let $1 \le i < j \le 2n$. Set
\begin{align*}
\beta_{i j} =
\begin{cases}
\alpha_{i} + \dots + \alpha_{j-1} + 2 \alpha_{j} + \dots +
2 \alpha_{2n-2} + \alpha_{2n-1} + \alpha_{2n} & (j < 2n)\\
\alpha_{i} + \dots + \alpha_{2n-2} + \alpha_{2n} & (j = 2n),
\end{cases}
\end{align*}
and $Y_{i j} = Y_{\beta_{i j}}$.
For a sequence $1 \le i_1 < i_2 < \dots < i_{2p} \le 2n$,
we set
\begin{align*}
(i_1, i_2, \dots, i_{2p}) = \sum_{\sigma \in \hat{S}_{2p}}
(- q^{-1})^{l(\sigma)} Y_{i_{\sigma(1), i_\sigma(2)}} \dotsm
Y_{i_{\sigma(2p-1), i_\sigma(2p)}},
\end{align*}
where $\hat{S}_m = \{ \sigma \in S_m \; | \;
\sigma(2k-1) < \sigma(2k+1), \sigma(2k-1) < \sigma(2k) \ 
\textrm{ for all } k \}$.
Then we have $f_{q,p} = (j_1^p, j_2^p, \dots, j_{2n}^p)$,
where $j_k^p=2n-2p+k$
(see \cite{Kamita}).
We can easily show the following description of $f_{q,p}$ similar
to \eqref{equat}.
\begin{align*}
f_{q,p} = \sum_{k=2}^{2p} (-q)^{2-k}
Y_{j_1^p, j_k^p} \;
(j_2^p, \dots, \check{j_k^p}, \dots, j_{2p}^p)
= \sum_{k=2}^{2p} (-q)^{2-k}
Y_{j_1^p, j_k^p} \;
\mathrm{ad}(F_{j_{k-1}^p} \cdots F_{j_2^p}) f_{q,p-1}.&
\end{align*}
Note that $\beta_{j_1^p, j_k^p} \notin \Delta_{(p-1)}^{+}$.
Hence we can
use Lemma \ref{lem of calulate-1}.

By using the induction on $i,j$, we can show the following lemma.
\begin{lem}
We have
\begin{align*}
{^t}Y_{j_k^p j_{k'}^p}(\partial) f_{q,p} = (-q)^{4n-1-k-k'}
(j_1^p, \dots, \check{j_k^p}, \dots, \check{j_{k'}^p}, \dots, j_{2p}^p)
\end{align*}
for $1 \le k < k' \le 2p$.
\end{lem}

Similarly to the case of type $A$, we obtain the following.
\begin{prop}
Let $(L_I, \mathfrak{n}_I^{-})$ be a regular prehomogeneous vector
space of type $(D_{2n}, 2n)$.
We have
\begin{align*}
a_p(s) = q^{\frac{s(4p+s-5)}{2}} \prod_{j=1}^s [j+2p-2]_q.
\end{align*}
In particular $b_{q,p}(s) = q^{s+2p-2} [s+2p-1]_q \; b_{q,p-1}(s)$,
we have the quantum $b$-function
\begin{align*}
b_{q,n}(s) =
\prod_{p=1}^n q^{s+2p-2} [s+2p-1]_q.
\end{align*}
\end{prop}

\vspace{5mm}

Let $(L_I, \mathfrak{n}_I^{+})$ be
the regular prehomogeneous vector space of type $(B_n, 1)$.
We label the vertices of the Dynkin diagram as in
Figure \ref{regular}, then $d_{i_0} = d_1 = 2$.
There exist two non-open orbits on $\mathfrak{n}_I^{+}$.
Let $1 \le i \le 2n-1$. We set $Y_i = Y_{\beta_i}$, where
\begin{align*}
\beta_i =
\begin{cases}
\alpha_1 + \dots + \alpha_i & (1 \le i \le n)\\
\alpha_1 + \dots + \alpha_{2n-i} + 2 \alpha_{2n-i+1} + \dots
+ 2 \alpha_n & (n+1 \le i \le 2n-1).
\end{cases}
\end{align*}
We have
\begin{align*}
&f_{q,1} = Y_1 = F_1,\\
&f_{q,2} = \sum_{i=1}^{n-1} (-q_{i_0})^{i+1-n} Y_{n+i} Y_{n-i} +
(q+q^{-1})^{-2} q^{-1} (-q_{i_0})^{1-n} Y_n^2
\end{align*}
(see \cite{Kamita}).
Note that $\beta_{i} \notin \Delta^{+}_{(1)}$ if $i \neq 1$.
On the other hand we have the following.
\begin{lem}
\begin{align*}
{^t}Y_i(\partial) f_{q, 2} =
\begin{cases}
(q+q^{-1})^{-1} (-q_{i_0})^{i-1} Y_{2n-i} & (1 \le i \le n)\\
-(q + q^{-1})^{-1}(-q_{i_0})^{i-2} Y_{2n-i} & (n+1 \le i \le 2n-1).
\end{cases}
\end{align*}
\end{lem}
Similarly to the case of type $A$, we obtain the following.
\begin{prop}
Let $(L_I, \mathfrak{n}_I^{-})$ be a regular prehomogeneous vector
space of type $(B_n, 1)$.
We have
\begin{align*}
a_2(s) = (q + q^{-1})^{-s}
q_{i_0}^{\frac{s(s+2n-4)}{2}} \prod_{i=1}^{s}
\left[ i + \frac{2n-3}{2} \right] _{q_{i_0}}.
\end{align*}
In particular we have the quantum $b$-function
\begin{align*}
b_{q,2}(s) = (q+q^{-1})^{-2} q_{i_0}^{s} \; [s+1]_{q_{i_0}}\; 
q_{i_0}^{s + \frac{2n-3}{2}} \;
\left[ s + \frac{2n-1}{2} \right] _{q_{i_0}}.
\end{align*}
\end{prop}

\vspace{5mm}

Let $(L_I, \mathfrak{n}_I^{+})$ be
the regular prehomogeneous vector space of type $(D_n, 1)$.
We label the vertices of the Dynkin diagram as in
Figure \ref{regular}, then $d_{i_0} = d_1 = 1$.
There exist two non-open orbits on $\mathfrak{n}_I^{+}$.
Let $1 \le i \le 2n-2$. We set $Y_i = Y_{\beta_i}$, where
\begin{align*}
\beta_i =
\begin{cases}
\alpha_1 + \dots + \alpha_i & (1 \le i \le n-1)\\
\alpha_1 + \dots + \alpha_{n-2} + \alpha_{n} & (i = n)\\
\alpha_1 + \dots + \alpha_{2n-i} + 2 \alpha_{2n-i+1} + \dots
+ 2 \alpha_{n-2} + \alpha_{n-1} + \alpha_n & (n + 1 \le i \le 2n-2).
\end{cases}
\end{align*}
Then we have
$f_{q,1} = Y_1 = F_1$, and
$f_{q,2} = \sum_{i=1}^{n-1} (-q)^{i+1-n} Y_{n+i-1} Y_{n-i}$
(see \cite{Kamita}).
We have the following results similar to those of type
$(B_n, 1)$.
\begin{lem}
\begin{align*}
{^t}Y_i(\partial) f_{q, 2} =
\begin{cases}
(-q)^{i-1} Y_{2n-1-i} & (1 \le i \le n-1)\\
(-q)^{i-2} Y_{2n-1-i} & (n \le i \le 2n-2).
\end{cases}
\end{align*}
\end{lem}

\begin{prop}
Let $(L_I, \mathfrak{n}_I^{-})$ be a regular prehomogeneous vector
space of type $(D_{n}, 1)$.
We have
\begin{align*}
a_2(s) = q^{\frac{s(s+2n-5)}{2}} \prod_{i=1}^{s} [i + n-2]_q.
\end{align*}
In particular we have the quantum $b$-function
\begin{align*}
b_{q,2}(s) = q^{s} [s+1]_{q} \; q^{s + n-2} \; [s + n-1]_q.
\end{align*}
\end{prop}

\vspace{5mm}

Let $(L_I, \mathfrak{n}_I^{+})$ be
the regular prehomogeneous vector space of type $(E_7, 1)$.
We label the vertices of the Dynkin diagram as in
Figure \ref{regular}, then $d_{i_0} = d_1 = 1$.
There exist three non-open orbits on $\mathfrak{n}_I^{+}$.

For $1 \le j \le 27$, we denote by $Y_j$ and $\psi_j$ the generators of
irreducible $U_q(\mathfrak{l}_I)$-modules
$V_q(\lambda_1)$ and $V_q(\lambda_2)$ respectively
(see \cite{Morita} for the explicit descriptions of $Y_j$ and $\psi_j$).
Note that $Y_j = Y_{\beta_j}$ for some $\beta_j
\in \Delta^{+} \setminus \Delta_I$,
$U_q(\mathfrak{n}_{(2)}^{-}) = \langle Y_1, \dots, Y_{10} \rangle$,
and $\psi_{27} = f_{q,2}$.
Now $(L_{(2)}, \mathfrak{n}_{(2)}^+)$ is of type
$(D_6, 1)$, hence we have
$b_{q,2}(s) = q^{s} [s+1]_{q} \; q^{s + 4} \; [s + 5]_q$.

The $q$-analogue $f_{q,3}$ of the basic relative invariant is given by
\begin{align*}
f_{q,3} &= \sum_{j=1}^{27} (-q)^{|\beta_j|-1} Y_j \; \psi_j\\
&= (1+q^8+q^{16}) Y_{27} \; \psi_{27} +
\frac{q^{-10}+q^{-8}-q^{-4}+1+q^{2}}{1+q^{2}}
\sum_{j=11}^{26} (-q)^{|\beta_{j}|-1} Y_{j} \; \psi_{j},
\end{align*}
where $|\beta| = \sum_{i=1}^7 m_i$ for
$\beta = \sum_{i = 1}^7 m_i \alpha_i$.
\begin{lem}
For $1 \le j \le 27$ we have
${^t}Y_j(\partial) f_{q,3} =
(1 + q^8 + q^{16}) (-q)^{|\beta_j|-1} \psi_j$.
\end{lem}
Then we have the following.
\begin{prop}
Let $(L_I, \mathfrak{n}_I^{-})$ be a regular prehomogeneous vector
space of type $(E_7, 1)$.
We have
\begin{align*}
a_3(s) = (1 + q^8 + q^{16})^{2s} q^{\frac{s(s+15)}{2}}
\prod_{i = 1}^{s} [i+8]_q.
\end{align*}
Therefore we have the quantum $b$-function
\begin{align*}
b_{q,3}(s) &=
(1 + q^8 + q^{16})^{2}\; q^{s+8} \; [s+9]_q \; b_{q,2}(s)\\
&= (1 + q^8 + q^{16})^{2}\; q^{s} [s+1]_{q} \; q^{s + 4} [s + 5]_q
\; q^{s+8} [s+9]_q.
\end{align*}
\end{prop}

\vspace{5mm}

Finally, we assume that $(L_I, \mathfrak{n}_I^{+})$ is
the regular prehomogeneous vector space of type $(C_n, n)$.
We label the vertices of the Dynkin diagram as in
Figure \ref{regular}, then $d_{i_0} = d_n = 2$.
There exist $n$ non-open orbits on $\mathfrak{n}_I^{+}$.
Let $1 \le i \le j \le n$.
We set $\beta_{i j} = \alpha_{i} + \dots + \alpha_{j-1} +
2 \alpha_{j} + \dots + 2 \alpha_{n-1} + \alpha_n$ and
$Y_{i j} = c_{i j} Y_{\beta_{i j}}$, where $c_{i j} = q + q^{-1}$
if $i = j$ and $1$ if $i \neq j$.
For $i < j$ we define $Y_{j i}$ by $Y_{j i} = q^{-2} Y_{i j}$.
Then we can write for $1 \le p \le n$
\begin{align*}
f_{q,p} = \sum_{\sigma \in S_p} (-q)^{-l(\sigma)}
Y_{i_1^{p}, i_{\sigma(1)}^{p}} \dotsm
Y_{i_p^{p}, i_{\sigma(p)}^{p}},
\end{align*}
where $i_k^{p} = n + k - p$ (see \cite{Kamita}).
\begin{lem}
\begin{align*}
f_{q,p} = Y_{i_1^{p},i_1^{p}} f_{q,p-1} + \sum_{k = 2}^p
\frac{(-q)^{1-k}}{q+q^{-1}} Y_{i_k^{p},i_1^{p}}
\ \mathrm{ad}(F_{i_{k-1}^{p}} \cdots F_{i_2^p}
F_{i_1^{p}})f_{q, p-1}.
\end{align*}
\end{lem}
\begin{proof}
We denote the right handed side of the statement by $g_p$.
It is easy to show that
the coefficient of $Y_{i_1^{p},i_1^{p}} \dotsm Y_{i_p^{p},i_p^{p}}$
in $f_{q,p}$ is equal to that in $g_p$.
Moreover the weight of $f_{q, p}$ is equal to that of $g_p$.
Hence it is sufficient to show that $g_p$ is the highest weight vector.
Since $(L_{(p)}, \mathfrak{n}_{(p)}^+) \simeq (C_p, p)$, we have only
to show the statement in the case where $p=n$.
We can easily show that $\mathrm{ad}(E_j) g_n = 0$ for $2 \le j \le n-1$.

Let us show $\mathrm{ad}(E_1) g_n = 0$.
For $2 \le j \le n$ we define $\varphi_j$ by
$\varphi_2 = \mathrm{ad}(E_1) g_n$ and
$\varphi_{j+1} = \mathrm{ad}(E_j) \varphi_j$.
We denote the weight of $\varphi_j$ by $\mu_j$.
Then we have $\mu_j \in - \alpha_1 + \sum_{i \neq 1}
\mathbb{Z}_{\le 0} \alpha_i$.
It is easy to show that
$\mathrm{ad}(E_k) \varphi_j = 0$ for any $k \neq j$.
In particular $\mathrm{ad}(E_k) \varphi_n = 0$ for any $k \in I$.

On the other hand we have the irreducible decomposition
\begin{align*}
U_q(\mathfrak{n}_I^{-}) = \bigoplus_{\mu \in \sum_{j=1}^{n}
\mathbb{Z}_{\ge 0} \lambda_j} V_q(\mu),
\end{align*}
and if $\mu \in \sum_{j=1}^{n} \mathbb{Z}_{\ge 0} \lambda_j$, then
$\mu \in 2 \mathbb{Z}_{\le 0} \alpha_1 + \sum_{i \neq 1}
\mathbb{Z}_{\le 0} \alpha_i$.
Hence $\mu_n \notin \sum_{j=1}^{n} \mathbb{Z}_{\ge 0} \lambda_j$,
and we have $\varphi_n = 0$.
We obtain $\varphi_j = 0$ for any $j$ by the induction.
\end{proof}
Note that $Y_{i_k^{p} i_{1}^{p}} \notin U_q(\mathfrak{n}_{(p-1)}^{-})$
for $1 \le k \le p$. Hence we can
use Lemma \ref{lem of calulate-1}.

We can prove the following lemma.
\begin{lem}
\begin{align*}
{^t}Y_{i_{1}^{p} i_k^{p}}(\partial) f_{q,p} =
\begin{cases}
(-q)^{2p-2} (q + q^{-1}) f_{q, p-1} & (k = 1)\\
-(-q)^{2p-k} \mathrm{ad} (F_{i_{k-1}^{p}} \cdots
F_{i_1^{p}})(f_{q, p-1}) & (k \le 2).
\end{cases}
\end{align*}
\end{lem}

Similarly to the case of type $A$, we obtain the following.
\begin{prop}
Let $(L_I, \mathfrak{n}_I^{-})$ be a regular prehomogeneous vector
space of type $(C_{n}, n)$.
We have
\begin{align*}
a_p(s) = (q + q^{-1})^s \; q_{i_0}^{\frac{s(s+p-2)}{2}}
\prod_{i=1}^{s} \left[ i + \frac{p-1}{2} \right]_{q_{i_0}}.
\end{align*}
In particular $b_{q,p}(s) = (q + q^{-1})\; q_{i_0}^{s + \frac{p-1}{2}}
\left[ s + \frac{p+1}{2} \right]_{q_{i_0}} \; b_{q,p-1}(s)$,
and we have the quantum $b$-function
\begin{align*}
b_{q,n}(s) = (q + q^{-1})^{n} \prod_{p=1}^{n}
q_{i_0}^{s + \frac{p-1}{2}} \left[ s + \frac{p+1}{2} \right]_{q_{i_0}}.
\end{align*}
\end{prop}

\vspace{20mm}
\begin{flushright}
\begin{tabular}{c}
\textit{Department of Mathematics}\\
\textit{Faculty of Science}\\
\textit{Hiroshima University}\\
\textit{Higashi-Hiroshima, 739-8526, Japan}
\end{tabular}
\end{flushright}
\end{document}